\documentclass[11pt,a4paper]{amsart}

\usepackage{amssymb,latexsym}
\usepackage{graphics}

\theoremstyle{plain}
\newtheorem{lemma}{Lemma}[section]
\newtheorem{proposition}[lemma]{Proposition}
\newtheorem{theorem}[lemma]{Theorem}
\newtheorem{corollary}[lemma]{Corollary}

\theoremstyle{remark}
\newtheorem{remark}[lemma]{Remark}

\theoremstyle{definition}

\numberwithin{equation}{section}

\newcommand{\R}{\mathbb{R}}

\newcommand{\Z}{\mathbb{Z}}

\newcommand{\TV}{\text{\rm Tot.Var.}}
\newcommand{\BV}{\text{\rm BV}}
\newcommand{\hyp}{\text{\rm Hyp}}
\newcommand{\M}{\mathbb{M}}

\begin{document}

\title[Stability of Riemann semigroup]{On the Stability of
the Standard Riemann Semigroup} 

\author{Stefano Bianchini}
\address{S.I.S.S.A. - I.S.A.S., via Beirut 2-4, 34014 TRIESTE (ITALY)}
\email{bianchin@sissa.it}
\urladdr{http://www.sissa.it/\textasciitilde bianchin/}

\author{Rinaldo M.~Colombo}
\address{Department of Mathematics, University of Brescia, Via Valotti
9, 25133 BRESCIA (ITALY)}
\email{rinaldo@ing.unibs.it}

\date{July 6th, 2000}
\thanks{We thank Alberto Bressan for useful discussions.}

\begin{abstract}
We consider the dependence of the entropic solution of a hyperbolic
system of conservation laws 
\[
\left\{ \begin{array}{c}
u_t + f(u)_x = 0 \\
u(0,\cdot) = u_0
\end{array} \right.
\]
on the flux function $f$. We prove that the solution in Lipschitz
continuous w.r.t.~the $C^0$ norm of the derivative of the perturbation
of $f$. We apply this result to prove the convergence of the solution
of the relativistic Euler equation to the classical limit.
\end{abstract}

\maketitle

% \vskip 1cm
\centerline{\small{Max-Planck-Institut f\"ur Mathematik in den Naturwissenschaften}}
\centerline{\small{Preprint n$^\circ$44 2000}}
\vskip 5mm

\section{Introduction}\label{S:SIintro}

Under suitable assumptions on the function $f\colon \Omega \mapsto
\R^n$ (with $\Omega \subseteq \R^n$), the system
\begin{equation}\label{E:Fsys}
u_t + \left[ f(u) \right]_x =0
\end{equation}
generates a {\it Standard Riemann Semigroup} (SRS) $S \colon \left[0,+\infty
\right[ \times \mathcal{D} \mapsto \mathcal{D}$, see~\cite{Br3}. Aim of
this paper is to investigate the dependence of $S$ upon the flow
function $f$.

Several papers in the current literature are concerned with the
existence of an SRS, see for example~\cite{BCP} and the references
in~\cite{Br3}. On the contrary, in the present paper the existence of an
SRS is assumed as a starting point and the focus is on the
correspondence $f \mapsto S$. In fact, the results in this
paper imply that the SRS $S$ is a Lipschitzean function of the flow $f$, with respect to the $C^0$ norm of $Df$. An immediate consequence is the following. Assume that $f$ depends on
the parameters $(p_1, \ldots, p_m)$ that may vary in a compact subset of
$\R^m$. Given a continuous functional $J$ defined on the solution $u$ at
time $t$ to the Cauchy problem for~\eqref{E:Fsys}, the present result
ensures the continuity of the map $(p_1,\ldots,p_m) \mapsto J(u(t))$
hence, by Weierstrass Theorem, the optimization problem admits a
solution.

For the sake of completeness, we only recall here that the existence of
the SRS for the $n\times n$ system \eqref{E:Fsys} was first proved in~\cite{BCP}.
The main assumptions there are that $Df$ is strictly hyperbolic with
every characteristic field either linearly degenerate or genuinely
nonlinear and that the initial data has sufficiently small total
variation. More recently, the existence of the SRS was extended also to
the non genuinely nonlinear setting in the $2\times2$ case, see~\cite{AM}.

Below we shall restrict our attention to {\it standard} solutions to
Riemann problems and, hence, to general Cauchy problems for Conservation
Laws. Here, by standard solutions we refer to those introduced by
Lax~\cite{Lax} and then generalized by Liu~\cite{Liu}. Various
extensions of the present work to other types of solvers are
straightforward.

The present paper is organized as follows. In the next section we state
the main results. The following Sections~\ref{S:SRel} and~\ref{S:Sstabinfty} are
devoted to two applications: the classical -- relativistic limit of
Euler equations and scalar conservation laws with $L^\infty$ initial
data. The proofs are given in Section~\ref{S:Stangvect}.

\section{Notation and Main Results}\label{S:SMain}

Consider the following hyperbolic system of conservation laws in one
space dimension
\begin{equation}\label{E:Fsysf}
u_t + \left[ f(u) \right]_x = 0 
\end{equation}
where $f\colon \Omega \mapsto \R^n$ is in $\hyp(\Omega)$, i.e.~$f$ is a
sufficiently smooth function that generate a SRS $S^f \colon
\left[0,+\infty \right[ \times \mathcal{D}^f \mapsto \mathcal{D}^f$. Recall that by SRS
generated by $f$ (see~\cite{Br1}) we mean a map $S^f \colon \left[0,+\infty
\right[ \times \mathcal{D}^f \mapsto \mathcal{D}^f$ with the following
properties:
\begin{itemize}
\item[(i)] $S^f$ is a semigroup: $S^f_0 = {\rm Id}$
and $S^f_t \circ S^f_s = S^f_{t+s}$;
\item[(ii)] $S^f$ is Lipschitz continuous: there exists a positive
$L_f$ such that for all positive $t,s$ and for all $u,w \in \mathcal{D}^f$,
$\| S^f_t u - S^f_s w \|_{L^1} \leq L_f \cdot \left(
|t-s| +\|u-w\|_{L^1} \right)$;
\item[(iii)] if $u$ is piecewise constant, then for $t$ small,
$S^f_t u$ coincides with the glueing of standard solutions to Riemann
problems.
\end{itemize}

\noindent For all $u \in \mathcal{D}^f$, it is well known (see~\cite{Br1}) that
the map $t \mapsto S^f_t u$ is a weak entropic solution to \eqref{E:Fsysf}.

Given $f\in\hyp(\Omega)$, let $\mathcal{R}(\mathcal{D}^f)$ be the set of all
piecewise constant functions in $\mathcal{D}^f$ having a single jump at the
origin. In other words, $\mathcal{R}(\mathcal{D}^f)$ is the set of initial
data to the Riemann problems
\begin{equation}\label{E:FRP}
\left\{ \begin{array}{c}
u_t + f(u)_x = 0 \\
u(0,x) = \begin{cases}	u^-	& \text{if } x<0	\\
			u^+	& \text{if } x>0	
\end{cases}	
\end{array}
\right.
\end{equation}
Below, by {\it solution} to \eqref{E:FRP} we always refer to the standard Lax
(see~\cite{Lax}) self--similar entropic solutions.

Let $f,g \in \hyp(\Omega)$ with
\begin{equation}\label{E:FIncl}
\mathcal{D}^g \subseteq \mathcal{D}^f
%\feqno{\FIncl}
\end{equation}
and define the {\it ``distance''} between $f$ and $g$ as (cfr.~\cite{Bi})
\begin{equation}\label{E:Fnormdist}
{\hat d}(f,g) =
\sup_{u \in \mathcal{R}(\mathcal{D}^g)}
\frac{1}{|u^+ - u^-|} \cdot \bigl\| S^f_1 u - S^g_1 u \bigr\|_{L^1}  
%\feqno{\Fnormdist}
\end{equation}
The distance ${\hat d}(f,g)$ is well defined due to \eqref{E:FIncl}. The main result of the
present paper is the following theorem.

\begin{theorem}\label{T:Tmain}
Let $f \in \hyp(\Omega)$. Then, for all $g \in
\hyp(\Omega)$ with $\mathcal{D}^g \subseteq \mathcal{D}^f$ and for all $u \in
\mathcal{D}^g$
\begin{equation}\label{E:Fmain}
\bigl\| S^f_t u-S^g_t u \bigr\|_{L^1}
\leq L_f
\cdot {\hat d} (f,g)
\cdot \int_0^t \TV \left( S^g_t u \right) \, dt \ .
\end{equation}
\end{theorem}

Recall that $L_f$ is the Lipschitz constant of the semigroup $S^f$,
see~(ii) above. The proof of Theorem~\ref{T:Tmain} is deferred to
Section~\ref{S:Stangvect}.

Remark that ${\hat d}$ generalizes the analogous quantity ${\hat
d}_{\text{lin}}$ defined in~\cite{Bi} with reference to the linear case.
Let $\M_d^{n\times n}$ denote the set of $n \times n$ diagonalizable
matrices with real eigenvalues. Note that $\M_d^{n\times n} \subseteq
\hyp(\R^n)$. Fix a $v \in \R^n$, $v\neq 0$. Denote by $A^t \star
v$ the solution evaluated at time $t$ of the linear system
\[
\left\{
\begin{array}{c}
u_t + A u_x = 0 				\\
u(0,x) = v \cdot \chi_{[0,+\infty)}(x)
\end{array}
\right.
\]
(here, $\chi_{I}$ is the characteristic function of the interval $I$).
Theorem~2.3 in~\cite{Bi} shows that
\begin{equation}\label{E:Fmatrixdist}
{\hat d}_{\text{lin}} (A,B)
\doteq
\sup_{v \colon |v|=1} \bigl\| A^1 \star v - B^1 \star v \bigr\|_{L^1},
\end{equation}
is a distance on $\M^{n\times n}_d$ such that for all $A,B \in
\M^{n\times n}$
\begin{equation}\label{E:Fparagi}
\| B-A \| \leq {\hat d}_{\text{lin}} (A,B) \ ,
\end{equation}
$\| \, \cdot \, \|$ being the usual operator norm. Moreover,
$\left(\M^{n\times n}_d, {\hat d}_{\text{lin}} \right)$ is a complete metric space.
Clearly, if $f$ and $g$ are linear, then ${\hat d}_{\text{lin}} (f,g) =
{\hat d}(f,g)$.

Furthermore, ${\hat d}$ is related to ${\hat d}_{\text{lin}}$ computed
on the derivatives of the flow functions, as shown by the following
Proposition.

\begin{proposition}\label{P:Pgeneral}
Let $f,g \in \hyp(\Omega)$ with $\mathcal{D}^g \subseteq \mathcal{D}^f$. Then
\begin{equation}\label{E:Fparagii}
{\hat d} (f,g)
\geq
\sup_{u \in \Omega} \ {\hat d}_{\text{{\rm lin}}} \! \left( Df(u), Dg(u) \right)  \ .
\end{equation}
\end{proposition}

Thus, ${\hat d}(f,g)$ seems stronger than the $C^0$ distance between
$Df$ and $Dg$, in the sense of \eqref{E:Fparagi}. Nonetheless,
Corollary~\ref{C:Ci}
below shows that once the flow functions and the domains $\mathcal{D}^f$,
$\mathcal{D}^g$ are fixed, i.e. the total variation of the solutions
$S^f_t u$, $S^g_t u$ are uniformly bounded, we can estimate the
r.h.s. in \eqref{E:Fmain} by means of $\| Df - Dg \|_{C^0}$. 

Theorem~\ref{T:Tmain} shows that the key point in the stability of the SRS
w.r.t.~the flow function lies in the dependence of only the solution to
Riemann problems upon the flow function. From the more abstract point of
view of quasidifferential equations in metric spaces (see~\cite{Br2,P})
this is equivalent to relate the distance between semigroups to the
distance between the vector fields generated by the semigroups.

As in~\cite{Br2} (see also~\cite{P}), in a metric space $(E,d)$ define an
equivalence relation on all the Lipschitz curves $\gamma \colon [0,1]
\mapsto E$ exiting a fixed point $u$ as
\begin{equation}\label{E:Feq}
\gamma \sim \gamma'
\quad \text{ iff } \quad
\lim_{\theta \mapsto 0^+} \frac{d \left( \gamma(\theta),\gamma'(\theta) \right)}{\theta} = 0.
\end{equation}
The quotient space $T_u$ so obtained is naturally equipped with the
metric
\begin{equation}\label{E:Fdist}
{\hat d}(v_1, v_2) \doteq
\limsup_{\theta \to 0^+} \frac{d \left( \gamma_1(\theta), \gamma_2(\theta) \right)}{\theta},
\end{equation}
where $\gamma_i$ is a representative of the equivalence class $v_i$.
By \eqref{E:Feq}, $\hat d$ does not depend on the particular representatives
chosen. A map $v\colon E \mapsto \bigcup_{u\in E}T_u$ is a vector field,
provided $v(u) \in T_u$ for all $u$.

Let $S \colon E \times [0,+\infty[ \mapsto E$ be a Lipschitzean
semigroup, i.e.~$S$ satisfies
\[
\begin{array}{c}
S_0 = \mathbb{I}_E		\\
S_s \circ S_t = S_{s+t}	
\end{array}
\quad \text{ and } \quad
\begin{array}{c}
\exists\, L>0 \text{ such that}				\\
\forall\, t,s>0 \text{ and } \forall\, u,w \in\mathcal{D}	\\
d \left( S_t u, S_s w \right)
\leq
L \cdot \left( |t - s| + d(u, w) \right) \ .
\end{array}
\]
Then, $S$ naturally defines a vector field $v_S$ on $E$ by
\begin{equation}\label{E:FvdaS}
v_S(u)
\quad \text{ is the equivalence class of the orbit } \quad
\theta \mapsto S_\theta u \ .
\end{equation}
Theorem~\ref{T:Tmain} has a natural abstract counterpart, namely

\begin{proposition}\label{P:Pcomp}
Let $S$, $S'$ be two Lipschitz
semigroups on $E$ generating the vector fields $v$ and $v'$, respectively.
Denote with $L$ the Lipschitz constant of, say, $S$. Then 
\begin{equation}\label{E:Ferror}
d \left( S_t u, S'_t u \right)
\leq
L \cdot \int_0^t {\hat d} \left( v(S'_t u), v'(S'_t u) \right) \, dt \ .
\end{equation}
\end{proposition}

The above proposition is an immediate corollary of the following widely
used (see~\cite{AM,Br1,Br2,BC,BCP} and the references
in~\cite{Br3}) error estimate:

\begin{lemma}\label{L:Lerrest}
Given a Lipschitz semigroup $S \colon E \times
[0,+\infty[ \mapsto E$ with Lipschitz constant $L$, for every Lipschitz
continuous map $w \colon [0,T] \mapsto E$ one has
\begin{equation}\label{E:Ferrest}
d \left( w(T), S_t w(0) \right)
\leq
L \cdot \int_0^T 
\liminf_{h \mapsto 0^+} \frac{d \left( w(t+h), S_h w(t) \right)}{h}
\ dt\ .
\end{equation}
\end{lemma}

\noindent The proof of Theorem~\ref{T:Tmain} consists of the following steps.
\begin{itemize}
\item[(1)] Find an explicit definition of the vector field $v_S$
generated by the SRS $S$.
\item[(2)] Compute the r.h.s.~in \eqref{E:Ferror}.
\item[(3)] Apply Lemma~\ref{L:Lerrest}.
\end{itemize}

\noindent Note that this procedure requires the mere existence of the
SRS. In several cases (see for instance~\cite{AM,BB,BC}) the
existence of the SRS is achieved through the construction of a sequence
$S^n$ of uniformly Lipschitzean approximate semigroups defined on
piecewise constant functions. In all these cases, the vector field $v_n$
generated by $S^n$ on the set of piecewise constant functions simply
consists in the gluing of solutions to Riemann problems. Thus, under
the further assumption that such an approximating sequence $S^n$ exists,
the step~(1) above could be avoided.

The quantity ${\hat d}$ in \eqref{E:Fnormdist} is thus the natural tool to
estimate the dependence of the SRS upon the flow function (note also
that no constant is involved in \eqref{E:Fmain}). However, in view of possible
applications of Theorem~\ref{T:Tmain}, we provide an estimate of the
r.h.s.~in \eqref{E:Fmain} in terms of handier quantities.

\begin{corollary}\label{C:Ci}
Let $f \in \hyp(\Omega)$ and assume that
\[
\mathcal{D}^f \subseteq \left\{ u \in L^{1}(\R,K) \colon \TV(u) \leq M \right\}
\]
for suitable positive $M$ and compact $K \subseteq \R^n$. Then,
there exists a positive constant $C$ such that for all $g \in
\hyp(\Omega)$ with $\mathcal{D}^g \subseteq \mathcal{D}^f$ and for all $u \in
\mathcal{D}^g$
\begin{equation}\label{E:Fcoroll} 
\bigl\| S^f_t u-S^g_t u \bigr\|_{L^1}
\leq C
\cdot \bigl\| Df-Dg \bigr\|_{C^{0}(\Omega)}
\cdot t \ .
\end{equation}
\end{corollary}

The above is the counterpart of the well known estimate for solutions of
bounded variation in the scalar case given in~\cite{BP}. Note that,
differently from the linear case, the ``distance'' $\hat d(f,g)$ is
equivalent to $\|Df - Dg\|_{C^0}$ because the total variation of both
solutions $S^f u$ and $S^g u$ is fixed by the domain $\mathcal{D}^f$: thus
the case of example of Remark 3.3 in~\cite{Bi} is not valid here.

For scalar
equations and assuming that the flow functions $f$ and $g$ are strictly
convex, we are able to extend the estimate in~\cite{BP} to $L^\infty$
initial data.

\begin{theorem}\label{T:Tinfty}
Assume that the scalar flow functions $f$ and
$g$ are uniformly strictly convex in a compact interval $K$:
i.e.~$f''(u), g''(u) \geq \kappa > 0$ for all $u \in K$. Let $u$
(resp.~$w$) denote the solution to
\[
\left\{
\begin{array}{c}
u_t + \left[ f(u) \right]_x =0	\\
u(0,x) = u_o(x),		
\end{array}
\right.
\quad \text{ resp.} \quad
\left\{
\begin{array}{c}
w_t + \left[ g(w) \right]_x =0	\\
w(0,x) = u_o(x)			
\end{array}
\right.
\]
with the same initial data $u_o \in L^\infty(\R, K)$. Denote by $\hat
\lambda$ an upper bound for the characteristic speeds,
i.e.~${\hat\lambda} \geq \max_{u\in K} \left\{ |f'(u)|,
|g'(u)| \right\}$. Then
\begin{equation}\label{E:Fhp}
\int_a^b \bigl| u(t,x) - w(t,x) \bigr| \, dx
\leq
2
\cdot \text{\rm diam}(K)
\cdot t
\cdot \frac{b - a + 4 {\hat\lambda} t}{\kappa t}
\cdot \max_{u \in K} \bigl| f'(u) - g'(u) \bigr| \ .
\end{equation}
\end{theorem}

\section{Proof of the main results}\label{S:Stangvect}

This section is devoted to the proofs of Theorem~\ref{T:Tmain},
Proposition~\ref{P:Pgeneral} and Corollary~\ref{C:Ci}. The first step
consists in the explicit computation of the vector field $v(u)$
generated by a SRS in the sense of \eqref{E:FvdaS}. This procedure
follows~\cite{Br2}. Remark that, given an $f \in \hyp(\Omega)$, such
construction is accomplished for all $u \in \mathcal{D}^f$.

Fix $f \in \hyp(\Omega)$ and $u \in \mathcal{D}^f$. $Du$ stands for the
total variation of the weak derivative of $u$. Let $\hat\lambda$ be a
constant strictly greater than all the characteristic speeds induced by
$f$ on $\mathcal{D}^f$. Let $\xi \in \R$ be given. We denote
by $\omega$ the self similar solution of the Riemann problem
\[
\left\{
\begin{array}{c}
\omega_t + f(\omega)_x =0	\\
\omega(0,x) = \begin{cases}
u(\xi-)	& \text{if } x<0	\\
u(\xi+)	& \text{if } x>0	
\end{cases}
\end{array}
\right.
\]
and let $U_{\xi}^{\sharp}$ be the function
\begin{equation}\label{E:Fusharp}
U_{\xi}^{\sharp}(\theta, x)
\doteq
\begin{cases}
\omega(\theta, x - \xi) & \text{if } |x - \xi| \leq	\hat \lambda \theta	\\
u(x)			& \text{if } |x - \xi| >	\hat \lambda \theta	
\end{cases}
\end{equation}
Moreover, define $U_{\xi}^{\flat}$ as the
solution to the linear hyperbolic problem with constant coefficients
\begin{equation}\label{E:Fuflat}
\left\{
\begin{array}{c}
\omega_{\theta} + Df(\xi) \omega_x = 0	\\
\omega(0, x) = u(x).
\end{array}
\right.
\end{equation}

Given $\epsilon > 0$ by an $\epsilon$--covering of the real line we
mean a family 
\begin{equation}\label{E:Fecover}
\mathcal{F}
\doteq
\left\{ I_1,\dots,I_N,I_1',\dots,I_M' \right\}
\end{equation}
of open intervals which cover $\R$ such that:
\begin{itemize}
\item[1)] the intervals $I_{\alpha}$ are mutually disjoint and no point
$x \in \R$ lies inside more than two distinct intervals
$I_{\beta}'$;
\item[2)] for every $\alpha$ there exists $\xi_{\alpha} \in I_{\alpha}$
such that $Du(I_{\alpha} \setminus \{ \xi_{\alpha} \} ) < \epsilon / N$;
\item[3)] for every $\beta$, $Du(I_{\beta}') < \epsilon$.
\end{itemize}
\noindent Now, let $\epsilon_n$ be a sequence strictly decreasing to
$0$. If $\mathcal{F}_n = \{I_1, \dots, I_{N(n)}, I_1', \dots,I_{M(n)}' \}$
is an $\epsilon_n$--covering of $\R$, then also $\mathcal{F}_{n,\theta}
= \{I_{1,\theta}, \dots, I_{N(n),\theta}, I_{1,\theta}', \dots,
I_{M(n),\theta}' \}$ is an $\epsilon_n$--covering, where $I_{j,\theta}
\doteq I_j + \lambda (-\theta,\theta)$ and $\theta \in [0, \theta_n]$,
for $\theta_n$ sufficiently small and such that the sequence $\theta_n$
strictly decreases to $0$. Define now
\[
u_n^{\theta}(x)
\doteq
\begin{cases}
U_{\xi_{\alpha}}^{\sharp}(\theta,x) & \text{if } x \in I_{\alpha,\theta}	\\
U_{\xi_{\beta}}^{\flat}(\theta,x) & \text{if } x \in
I_{\beta,\theta}', x \notin \bigcup_{\alpha} I_{\alpha,\theta}, x
\notin \bigcup_{\beta' < \beta} I_{\beta',\theta}' 
\end{cases}
\]
We finally obtain the vector field $v_{S^f}$ letting $v_{S^f}(u)$ be the
equivalence class (w.r.t.~\eqref{E:Feq}) of the curve $\theta \mapsto u^\theta$,
where
\begin{equation}\label{E:Fuenne}
u^{\theta}
\doteq
s u_n^{\theta_n} + (1 - s) u_{n-1}^{\theta_{n-1}}
\quad \text{ if } \quad
\theta = s \theta_{n} + (1 - s) \theta_{n-1}, \ s \in [0,1] \ .
\end{equation}

In~\cite{Br2} it is shown that the trajectories of SRS $S^f$ generated
by the hyperbolic system of conservation laws $u_t + f(u)_x=0$ are the
solution of the quasilinear equation ${\dot u} = v_{S^f}(u)$, where the
vector field $v_{S^f}(u)$ is generated by the curve~\eqref{E:Fuenne}.

Let $f$, $g$ be as in Theorem~\ref{T:Tmain}, and denote with $u^{\theta}$,
$w^{\theta}$ the two curves generating the vector fields $v_{S^f}$ and
$v_{S^g}$ induced by the SRSs $S^f$ and $S^g$. In view of \eqref{E:Ferror}, we
now pass to compute the distance \eqref{E:Fdist} between $v_{S^f} (u)$ and
$v_{S^g} (u)$. A simple computation shows that by
Proposition~\ref{P:Pgeneral}
\begin{align}
\frac{d ( u^{\theta_n}, w^{\theta_n} )}{\theta_n}
& =
\frac{1}{\theta_n} \cdot \bigl\| w^{\theta_n} - u^{\theta_n}
\bigr\|_{L^1} \\
& \leq
\sum_{\alpha} Du(\xi_{\alpha}) \cdot {\hat d}_{\text{lin}} (f,g) +
\sum_{\beta} \TV(u,I_{\beta}') \cdot {\hat d} \left(Df(\xi_{\beta}), Dg(\xi_{\beta}) \right)
\notag \\
& \leq
\TV(u) \cdot {\hat d} (f,g) \ .
\notag 
\end{align}
As a consequence, if $s$ is as in \eqref{E:Fuenne}
\begin{align}
\frac{d(u^\theta,w^\theta)}{\theta}
&=
\frac{\bigl\| s\cdot (u^{\theta_n}_n - w^{\theta_n}_n) + (1-s) \cdot
(u^{\theta_{n-1}}_{n-1} - w^{\theta_{n-1}}_{n-1}) \bigr\|_{L^1}}{\theta}
\\
& \leq
\frac{s \cdot \theta_n + (1-s) \cdot \theta_{n-1}}{\theta} \cdot \TV(u) \cdot {\hat d} (f,g)
\notag \\
& =
\TV(u) \cdot {\hat d} (f,g)
\notag 
\end{align}
hence
\begin{equation}\label{E:Ffinal}
{\hat d} \left( v_{S^f} (u), v_{S^g} (u) \right)
=
\limsup_{n\to+\infty} \frac{d ( u^\theta, w^\theta )}{\theta}
\leq
\TV(u) \cdot {\hat d} (f,g) \ .
\end{equation}
By Proposition~\ref{P:Pcomp}, applying \eqref{E:Ffinal} to $S^g_tu$,
the proof of Theorem~\ref{T:Tmain} is completed.

In the next part of this section we give a proof of
Corollary~\ref{C:Ci}. We only need to prove that there exists a
constant $C$ such that $\hat d(f,g) \leq C \cdot \|Df - Dg\|_{C^0}$:
this means that for all Riemann problems \eqref{E:FRP} we have
\begin{equation}\label{E:cazzo1}
\bigl\| S^f_1 u - S^g_1 u \bigr\|_{L^1} \leq C \cdot \bigl\|Df - Dg \bigr\|_{C^0}
\cdot \bigl| u^+ - u^- \bigr|.
\end{equation}
We recall that $S^g_t u$ is a self similar solution, obtained by
piecing together centered rarefaction waves and jump
discontinuities. By the $L^1_{\text{loc}}$ dependence of $S^f$,
formula \eqref{E:cazzo1} is proved if we can verify it for Riemann
problems generating a single wave in the solution $S^g_t u$. We then
have to consider two cases.

If $S^g_t u$ is a centered rarefaction wave, then by the Lipschitz
continuity of the solution for $t > 0$, the functions $S^g_{t+h} u$
and $S^f_h \circ S^g_t u$ solves in the broad sense the quasilinear
versions of \eqref{E:Fsysf}:
\begin{equation}\label{E:cazzo2}
\bigl[S^f_h \circ S^g_t u \bigr]_h + Df \bigl[ S^f_h \circ S^g_t u
\bigr]_x = 0 \quad \text{and} \quad \bigl[S^g_{t+h} u \bigr]_h
+ Dg \bigl[ S^g_{t+h} u \bigr]_x = 0.
\end{equation}
Applying Lemma~\ref{L:Lerrest}, we obtain
\begin{align}\label{E:cazzo3}
\bigl\| S^f_1 u - S^g_1 u \bigr\|_{L^1} \leq&~ L_f \cdot \int_0^1
\liminf_{h \to 0} \frac{\bigl\| S^f_h  \circ S^g_t u - S^g_{t+h} u
\bigr\|_{L^1}}{h} dt \\
=&~ L_f \cdot \int_0^1 \bigl\| Df \bigl[S^g_t u\big]_x -
Dg \bigl[ S^g_{t} u \bigr]_x \bigr\|_{L^1} \\
\leq&~ C \cdot \bigl\|Df - Dg
\bigr\|_{C^0} \cdot |u^+ - u^-|, \notag
\end{align}
since in this case the $L^1$ limits as $h \to 0$ of $(S^f_h  \circ
S^g_t u)/h$ and $(S^g_{t+h} u)/h$ exist by \eqref{E:cazzo2} and are
equal to $(Df - Dg) [S^g_{t} u ]_x$. 

Now we consider the case in which the jump $u^-, u^+$ is solved by a
shock travelling with speed $\sigma$, where $\sigma$ is given by the
Rankine--Hugoniot condition
\begin{equation}\label{cazzo4}
g(u^+) - g(u^-) = \sigma \cdot \bigl( u^+ - u^- \bigr).
\end{equation}
To prove \eqref{E:cazzo1}, we approximate the
solution $S^g_t u$ with the Lipschitz continuous function $\tilde
u(t)$ defined as
\[
\tilde u(t,x) = u^- \cdot \chi_{(-\infty,0]} \bigl( x - \sigma t
\bigr) + \bigl( u^+ - u^- \bigr) \cdot \min \left( \frac{x - \sigma
t}{\delta}, 1 \right) \cdot \chi_{(0,+\infty)} \bigl(x - \sigma t
\bigr). \]
Roughly speaking, $\tilde u$ is obtained from $S^g_t u$ by substituting
the jump $u^-,u^+$ at $x - \sigma t$ with a linear function. Using the
Lipschitzeanity of $S^g$ we can write
\begin{align}\label{E:Ftrucco}
\bigl\| S_h^f \circ S^g_t u - S^g_{t+h} u \bigr\|_{L^1}
& \leq
\bigl\| S_h^f \circ S^g_t u - S^f_h \tilde u(t) \bigr\|_{L^1}
+
\bigl\| S^f_h \tilde u(t) - \tilde u(t+h) \bigr\|_{L^1}
\\
& \qquad +
\bigl\| \tilde u(t+h) - S_{t+h}^g u \bigr\|_{L^1}
\notag \\
&\leq
\frac{\delta}{2} \cdot (1 + L^f) \cdot \bigl| u^+ -
u^- \bigr| +
\bigl\| \tilde u(t+h) - S_h^f \tilde u(t)
\bigr\|_{L^1} \, . \notag
\end{align}
The last term above can be evaluated again by means of Lemma~\ref{L:Lerrest}:
\begin{align}
&\bigl\| \tilde u(t+h) - S_h^f \tilde u(t) \bigr\|_{L^1}
\leq L^f \cdot
\int_t^{t+h} \liminf_{\xi \to 0^+}
\frac{\bigl\| \tilde u(\eta + \xi) - S_\xi^f \tilde u(\eta)
\bigr\|_{L^1}}{\xi} \, d\eta
\\
& \qquad =
h \cdot L^f\cdot
\int_0^{\delta}
\Bigl|
	\left( \dot \sigma_{\alpha} \mathbb{I} -
		Dg\left( u_{\alpha}^- + (u_{\alpha}^+ - u_{\alpha}^-) \frac{y}{\delta} \right)
	\right)
	\left( u_{\alpha}^- + (u_{\alpha}^+ - u_{\alpha}^-) \frac{y}{\delta} \right)
\Bigr| \, dy
\notag \\
& \qquad \leq
h \cdot L^f
\cdot \sup_{x \in K} \bigl\| Df(x) - Dg(x) \bigr\|
\cdot \bigl| u_{\alpha}^+ - u_{\alpha}^- \bigr| \, , \notag
\end{align}
where we use $\dot x_{\alpha} \in [\sigma - \epsilon, \sigma +
\epsilon]$ and the relation
\[
\sigma \cdot (u_\alpha^+ - u_\alpha^-)
=
f(u_\alpha^+) - f(u_\alpha^-)
=
\int_0^1 Df \left( (1-s) u_\alpha^- + s u_\alpha^+ \right) ds \cdot (u_\alpha^+ - u_\alpha^-) \, .
\]
In fact, since $u_{\alpha,\delta}$ is Lipschitz continuous, we can use
Lebesgue's dominated convergence theorem as in the previous case.
Letting $\delta$ tend to $0$, we obtain finally
\begin{equation}\label{E:Ferrshock}
\frac{1}{h} \bigl\| S^f_h  \circ S^g_t u - S_{t+h}^g \bigr\|_{L^1}
\leq L_f \cdot \sup_{x \in K} \|Df(x) - Dg(x)\| \cdot \bigl|u_{\alpha}^+
- u_{\alpha}^- \bigr|.
\end{equation}
This concludes the proof of Corollary~\ref{C:Ci}: in fact an
application of Lemma~\ref{L:Lerrest} gives immediately~\eqref{E:cazzo1}.

To end this section, we prove Proposition~\ref{P:Pgeneral}.
Fix $u_o \in \Omega$ and $v \in \R^n$
with $|v|=1$. Assume that for all positive and sufficiently small
$h$ the function $u_h=u_o + h \cdot \chi_{[0,+\infty)}\cdot v$ is in
$\mathcal{D}^g$. Let
\[
{\tilde f}_h (u)= \frac{1}{h} \cdot \bigl( f(u) - f(u_o) \bigr)
\ , \qquad
{\tilde g}_h (u)= \frac{1}{h} \cdot \bigl( g(u) - g(u_o) \bigr)
\]
and note that using a simple rescaling, we can write
\[
\left\| S^{{\tilde f}_h}_1 \left( \frac{1}{h}(u_h-u_o) \right) -
	S^{{\tilde g}_h}_1 \left( \frac{1}{h}(u_h-u_o) \right) \right\|_{L^1}
=
\frac{1}{h} \cdot \bigl\| S^f_1 u_h - S^g_1 u_h \bigr\|_{L^1}
\]
(recall that $(1/h)\cdot(u_h-u_o) = \chi_{[0,+\infty)} \cdot v$ is
independent from $h$). Hence, passing to the limit $h \to 0$
\begin{align}\label{E:Fdistlimit}
\bigl\| Df(u_o)^1 \star v - Dg(u_o)^1 \star v \bigr\|_{L^1}
& =
\lim_{h\to 0} \frac{1}{h} \cdot \bigl\| S_1^f u_h - S_1^g u_h \bigr\|_{L^1}
\\
{\hat d}_{\text{{\rm lin}}} \left( Df(u_o), Dg(u_o) \right)
& =
\sup_{v}
\lim_{h\to 0} \frac{1}{h} \cdot  \bigl\| S_1^f u_h - S_1^g u_h \bigr\|_{L^1}
\notag
\end{align}
and the proof is completed.

\section{The Classical Limit of the Relativistic Euler Equations}\label{S:SRel}

In this section we apply Corollary~\ref{C:Ci} to the classical limit of the
relativistic Euler equations, generalizing what was obtained in~\cite{MU}.

The relativistic $p$--system, (see~\cite{ST,T}) is
\begin{equation}\label{E:FRelEul}
\left\{
\begin{array}{c}
{\displaystyle \left[ \rho + \left(\rho + \frac{1}{c^2} p\right)
\frac{(v/c)^2}{1-(v/c)^2} \right]_t + 
\left[ \left(\rho + \frac{1}{c^2} p\right) \frac{v}{1-(v/c)^2} \right]_x
=0 } \\[11pt]
{\displaystyle \left[ \left(\rho + \frac{1}{c^2} p\right) \frac{v}{1-(v/c)^2} \right]_t +
\left[ \left(\rho + \frac{1}{c^2} p\right) \frac{v^2}{1-(v/c)^2} +p\right]_x
=0 }
\end{array}
\right.
\end{equation}
Above, $\rho$ is the mass--energy density of the fluid, $v$ the
classical coordinate velocity, $p$ the pressure and $c$ the light speed.
We show below that as $c\to+\infty$, the problem \eqref{E:FRelEul} approaches to
its classical counterpart
\begin{equation}\label{E:FEul}
\left\{
\begin{array}{c}
\left[ \rho \right]_t + \left[ \rho v \right]_x =0	\\[4pt]
\left[ \rho v \right]_t + \left[ \rho v^2+p \right]=0 	
\end{array}
\right.
\end{equation}
in the sense that the SRS $S^c$ generated by \eqref{E:FRelEul} converges to the
SRS $S$ generated by \eqref{E:FEul} on a domain containing all physically
reasonable data. In particular, the total variation of the data need not
be small.

In~\eqref{E:FRelEul} a standard choice (see~\cite{ST} and the references
therein) for the pressure law is
\begin{equation}\label{E:FpLaw}
p= \sigma^2 \cdot \rho \ ,
\end{equation}
$\sigma$ being the sound speed. 

Fix a positive lower bound for the density $\rho_{\min}$ and for the
light speed $c_o$. Without any loss in generality, we may assume
$\sigma<c_o$.

Let $\mathcal{V}_M =\left\{ (\rho,\rho v) \in \BV \left( \R,
(\rho_{\min},+\infty) \times \R \right) \colon \TV(\rho) + \TV(\rho
v) \leq M\right\}$. In~\cite{CR} it is proved that for any $M>0$, \eqref{E:FRelEul}
generates a SRS $S^{c,M}$ defined on a domain $\mathcal{D}^{c,M}$
containing $\mathcal{V}_M$ and consisting of functions of total variation
bounded by, say, $M^3$ (provided $M$ is sufficiently large). Similarly,
\eqref{E:FEul} generates a SRS $S^{M}$ on a domain $\mathcal{D}^{M}$ containing
$\mathcal{V}_{M^3}$ and contained, say, in $\mathcal{V}_{M^9}$.

Thus, for all sufficiently large $c$ and $M$, there exist domains
$\mathcal{D}^{c,M}$, $\mathcal{D}^{M}$ such that 
\begin{equation}\label{E:FIncl1}
\mathcal{V}_M
\subset \mathcal{D}^{c,M}
\subset \mathcal{V}_{M^3}
\subset \mathcal{D}^{M}
\subset \mathcal{V}_{M^9}
\end{equation}
and, moreover, the problems \eqref{E:FRelEul} and \eqref{E:FEul} generate the SRSs
\[
S^{c,M} \colon \mathcal{D}^{c,M} \times (0,+\infty) \mapsto \mathcal{D}^{c,M}
\quad \text{ and } \quad 
S^{M} \colon \mathcal{D}^{M} \times (0,+\infty) \mapsto \mathcal{D}^{M}
\ .
\]

We are now ready to state and prove the following application of
Corollary~\ref{C:Ci}.

\begin{theorem}\label{T:TRel}
Fix a positive $\rho_{\min}$ and sufficiently large
$M$, $c_o$. Let $\mathcal{D}^{c,M}$ and $\mathcal{D}^{M}$ satisfy to \eqref{E:FIncl1}
for $c>c_o$. Then there exists a constant $C$ such that for all $c \geq
c_0$ and for all $u \in \mathcal{D}^{M,c}$
\begin{equation}\label{E:Fsemdit}
\bigl\| S^{c,M}_t u - S^{M}_t u \bigr\|_{L^1}
\leq
C \cdot \frac{1}{c^2} \cdot t.
\end{equation}
\end{theorem}

In particular, by \eqref{E:FIncl1}, the bound \eqref{E:Fsemdit} holds for all initial
data $u$ with $\TV(u) < M$.

\begin{proof}
Note that \eqref{E:FRelEul} and \eqref{E:FEul} in conservation form become,
respectively,
\begin{equation}\label{E:FDueSis}
\left\{
\begin{array}{c}
\rho_t + q_x 							=0	\\
{\displaystyle q_t + \left( \phi_c(\rho,q) \cdot \frac{q^2}{\rho}+ p
\right)_x=0 }
\end{array}
\right.
\quad \text{ and } \quad
\left\{
\begin{array}{c}
\rho_t + q_x 					=0	\\
{\displaystyle q_t + \left( \frac{q^2}{\rho} + p \right)_x	=0 }
\end{array}
\right.
\end{equation}
where
\[
\phi_c(\rho,q) =
1+ \frac{1}{c^2}
\cdot \left( 1-\frac{v^2(\rho,q)}{c^2} \right)
\cdot \frac{p}{{\displaystyle \rho+ \frac{v^2(\rho,q)}{c^2} \cdot
\frac{p}{c^2}}}
\ .
\]
Call $f_c$ and, respectively, $f$ the fluxes in the two
systems \eqref{E:FDueSis}. Then, the estimate $\| Df_c - Df \|_{C^0} \leq C
\cdot (1/c^2)$ follows from straightforward computations and completes
the proof.
\end{proof}

We remark that the particular pressure law \eqref{E:FpLaw} is necessary only to
ensure the existence of the SRS $S^{c,M}$ in the large. The above
procedure remains true under much milder assumptions on the equation of
state.
Note moreover that the rate of convergence $\mathcal{O}(1/c^2)$ is exactly the one
expected by the convergence of relativistic to classical mechanics.

It is of interest to mention that~\eqref{E:Fsemdit} proves also the uniform
Lipschitz continuity of all semigroups $S^{c,M}$ for all
sufficiently large $c$.

\begin{remark}

As is well known, linearizing~(\ref{E:FEul}) around $\rho=\rho_o$, $v=0$ at constant entropy leads to the wave equation. Corollary~\ref{C:Ci} ensures that the solutions to~(\ref{E:FEul}) converge to the linearized equation in $L^1$ over finite time intervals.

\end{remark}

\section{Stability of a scalar equation w.r.t.~flux}\label{S:Sstabinfty}

In this section we prove Theorem~\ref{T:Tinfty}. We consider two scalar equation,
\begin{equation}\label{E:Fequlo}
u_t + f(u)_x = 0,
\end{equation}
\begin{equation}\label{E:Fequ}
v_t + g(u)_x = 0,
\end{equation}
with the same initial condition: $u(0,x) = v(0,x) = u_0$. We assume
that $f$, $g$ are strictly convex $C^2$ functions, precisely there exists a
constant $\kappa$ such that
\begin{equation}\label{E:Fcd}
\min_{u \in K} \{|f''(u)|, |g''(u)| \} \geq \kappa \, ,
\end{equation}
where $K$ is a compact interval of $\R$ such that $u_0(x) \in K$ for
all $x \in \R$. We recall that by maximum principle also the entropic
solutions of~\eqref{E:Fequlo}, \eqref{E:Fequ} will satisfies the same
bound.

We recall that, by~\cite{Lax}, given a point $(t,x)$ we can consider the
set of characteristics $\xi(t,x)$ passing through $(t,x)$. If we
denote with $\xi^-(t,x)$ and $\xi^+(t,x)$ the minimal and maximal
backward characteristics, then either $\xi^-(t,x) = \xi^+(t,x)$ and
the solution $u$ is continuous in $(t,x)$, or we have an admissible shock and the
jump is exactly given by the (constant) values of $u$ on the
characteristics $\xi^-$, $\xi^+$. By condition~\eqref{E:Fcd} we have
that if at time $t$ two characteristics $\xi^+(t)$ and $\xi^-(t)$ meet,
then we have \begin{equation}\label{E:Fch}
\frac{d}{dt} (\xi^+ - \xi^-) \leq - \kappa (u^+ - u^-) \, .
\end{equation}

Suppose that $\xi^-(0) < \xi^+(0)$, and consider now an initial datum
$\tilde u_0$ defined as
\begin{equation}\label{E:Find}
\widetilde u_0(x) = \begin{cases}
u_0(x) & x \leq \xi^-(0) \\
u_0^- & \xi^-(0) < x \leq \Xi \\
u_0^+ & \Xi < x \leq \xi^+(0) \\
u_0(x) & x > \xi^+(0) \end{cases} 
\end{equation}
where $\Xi \in \R$ is chosen such that 
\begin{equation}\label{E:Fpos}
\Xi \doteq \frac{1}{u_0^- - u_0^+} \int_{\xi^-}^{\xi^+} u_0(x) dx +
\frac{u_0^- \xi^- - u_0^+ \xi^+}{u_0^- - u_0^+} \,.
\end{equation}
Using the conservation of mass it is easy to conclude that the
solution in unchanged at time $t$. In fact consider the triangle $\mathcal{T}$ whose 
vertices are $(0,\xi^-(0))$, $(0, \xi^+(0))$ and $(t,x)$. Since the
equation \eqref{E:Fequlo} can be written as
\[
\text{\rm div} \left( \begin{array}{c} u \\ f(u) \end{array} \right) = 0 \, ,
\]
and $u$ is constant along the lines $\xi(t)$, we have
\begin{align}
\int\!\!\!\int_\mathcal{T} (u_t + f(u)_x) dtdx =& - \int_{\xi^-}^{\xi^+}
u_0(x) dx + \int_{0}^t (f(u^+) - f'(u^+)u^+ ) dt + \int_{0}^t
(f'(u^-)u^- - f(u^-)) dt \notag \\
=& - \int_{\xi^-}^{\xi^+} u_0(x) dx + (f'(u^-)u^- - f'(u^+)u^+ + f(u^+)
- f(u^-)) t = 0. \notag
\end{align}
Using the relation
\[
x = \xi^- + f'(u^-) t = \Xi + \frac{f(u^+) - f(u^-)}{u^+ - u^-}t =
\xi^+ + f'(u^+) t \, ,
\]
we obtain \eqref{E:Fpos}. Using \eqref{E:Fpos} we can change the initial
of~\eqref{E:Fequlo} or~\eqref{E:Fequ} so that $u(t)$ or $v(t)$ are
unchanged. Since with this procedure we collect all the interactions at
time $0$, the total variation of $\tilde u_0$ has the same value of
$\TV(u(t))$.

To change the initial datum in such a way that both $u(t)$ and $v(t)$
are the same, consider now the test system 
\begin{equation}\label{E:Ftest}
w_t + \left(\frac{\kappa}{2} w^2 \right)_x = 0 \,.
\end{equation}
By \eqref{E:Fcd} and \eqref{E:Fch} if two characteristics meet is $u(t)$ and $v(t)$,
then they also meet in \eqref{E:Ftest}. Let us denote with $\tilde u_0$ the new
initial condition, obtained by the above procedure using equation
\eqref{E:Ftest}.

Consider now an interval $[a,b]$. By the definition of ${\hat\lambda}$ we have that the values of
$u(t)$ and $v(t)$ in $[a,b]$ depends only on $\tilde u_0$ in $[a -
\hat \lambda t, b + \hat \lambda t]$. Using the standard estimates we have 
\begin{equation}\label{E:Fint1}
\int_a^b |u(t,x) - v(t,x)| dx \leq t \cdot \max_{u \in K} |f'(u) -
g'(u)| \cdot \TV \Bigl( \widetilde u_0;[a - \widehat \lambda t, b +
\widehat \lambda t] \Bigr).
\end{equation}
To estimate the total variation, an easy computation gives
\begin{equation}\label{E:Ftv}
\TV \Bigl( \widetilde u_0;[a - \widehat \lambda t, b + \widehat
\lambda t] \Bigr)
\leq \TV \Bigl( w(t);[a - 2 \widehat \lambda t, b + 2 \widehat \lambda
t] \Bigr)
\leq 2 \cdot \text{diam}(K) \cdot \frac{b - a + 4 \widehat \lambda
t}{\kappa t} \,.
\end{equation}
Combining~\eqref{E:Fint1} and~\eqref{E:Ftv} we get
\begin{equation}\label{E:Ffin}
\int_a^b |u(t,x) - v(t,x)| dx \leq \max_{u \in K} \bigl| f'(u) - g'(u)
\bigr| \cdot 2 \cdot
\text{diam}(K) \cdot \frac{b - a + 4 \widehat \lambda t}{\kappa
}.
\end{equation}

\begin{remark}\label{E:Rexp}
When $t \to 0$ the integral does not converge
to $0$. This is clear since the initial datum is in $L^{\infty}$, and
then the semigroup is continuous but not Lipschitz continuous in time,
since the amount of interaction at $t = 0$ is infinite. Consider for
example the following two equations 
\[
u_t + \left( \frac{u^2}{2} \right)_x = 0,
\qquad \qquad
v_t + \left( - v + \frac{v^2}{2} \right)_x = 0 \ ,
\]
with the periodic initial condition
\[
u_{0,n}(x)
=
\begin{cases}
1	&	\text{if }x \in \bigl[k 2^{- n + 1}, k 2^{- n + 1} + 2^{ -
n} \bigr],k \in \Z
\\
-1	&	\text{otherwise}
\end{cases}
\]
At time $t_n = 2^{-n}$ the two solution are
\[
\begin{cases}
u_n (t_n,x) = 2^n (x - k 2^{-n + 1}) 
& \text{if } x \in \bigl( k 2^{- n + 1} - 2^{ - n}, k 2^{- n + 1} +
2^{ - n} \bigr]
\\
v_n (t_n,x) = 2^n (x - k 2^{-n + 1} - 2^{-n})
& \text{if } x \in \bigl( k 2^{- n + 1}, (k + 1) 2^{- n + 1} \bigr] 
\end{cases}
\]
so that
\[
\int_{0}^1 \bigl| v_{n}(t_n,x) - u_n(t_n,x) \bigr| \, dx = 1 \ .
\]
This depends on the fact that the modulus of continuity of the semigroup can be
arbitrarily large.

Note moreover that since the solutions $u$, $v$ are limits of wave
front tracking approximations, the continuous dependence of the
solution on the flux function $f$ can be stated also if $f$ non
convex. However in general one cannot prove any uniform continuous
dependence. 
\end{remark}

\medskip

\noindent\textbf{Acknowledgments.}  The present work was partly accomplished while the authors visited the Max Planck Institute in Leipzig supported by the European TMR Network \emph{``Hyperbolic Conservation Laws''} {ERBMRXCT960033}.

\end{document}